\documentclass[11pt]{article}
\usepackage{customtemplate}
\usepackage[margin=1.15in]{geometry}
\usepackage{amsmath,amssymb,amsthm}

\theoremstyle{plain}

\theoremstyle{definition}

\usepackage{enumitem}
\newcommand{\No}{\mathbb{N}_{\mathrm{o}}}          
\newcommand{\dodd}{d_{\mathrm{o}}}                 
\newcommand{\R}{\mathfrak{R}}
\newcommand{\odd}{\operatorname{od}}
\newcommand{\wo}{\omega_{\mathrm{o}}}              
\newcommand{\N}{\mathbb{N}}

\title{\textbf{A Divisor-Sum Analogue of the Collatz Map}}
\author{%
Ritesh Dwivedi\footnote{University of Allahabad, Prayagraj, India. Email: \href{mailto:riteshgrouptheory@gmail.com}{riteshgrouptheory@gmail.com}}
\and
Rohit Yadav\footnote{Indian Institute of Technology Jammu, Jammu, India 181221.
Email: \href{mailto:rohityadavau1998@gmail.com}{rohityadavau1998@gmail.com}}}
\date{}

\begin{document}
\maketitle

\begin{abstract}
Let $\sigma$ denote the sum-of-divisors function and let $\R$ send an odd integer $n$ to $\sigma(n)$ and an even integer $n$ to $n/2$. We conjecture that every orbit of $\R$ reaches $1$; this implies that there is no odd $2^{k}$-perfect number for any $k \ge 1$. We prove a pointwise descent estimate for the map $T$ induced by $\R$ on the odd integers, and show that a single application of $T$ divides almost every odd $n$ by $(\log n)^{2\log 2-\varepsilon}$, for every fixed $\varepsilon>0$. The estimate does not iterate, and we determine what is missing.
\end{abstract}



\section{Introduction}\label{sec:intro}

Among the recursively defined integer sequences studied in elementary number theory, those of Collatz type have received particular attention. The classical Collatz map halves an even integer and sends an odd integer $n$ to $3n+1$, and the conjecture that each of its orbits reaches $1$ has remained open for several decades; see, for example, \cite{Lagarias2010, Tao2022} and the relevant references therein. In this paper we consider a variant of this map in which the affine rule on the odd integers is replaced by an arithmetic one, namely the sum-of-divisors function.

Let $\sigma(n)=\sum_{d \mid n} d$. For a positive integer $n$, define
\begin{equation}\label{eq:R}
\R(n):=
\begin{cases}
\sigma(n), & \text{if $n$ is odd},\\[2pt]
n/2, & \text{if $n$ is even},
\end{cases}
\end{equation}
and set $\R^{0}(n):=n$ and $\R^{k+1}(n):=\R(\R^{k}(n))$ for $k \ge 0$. Since $\sigma(1)=1$, the integer $1$ is a fixed point of $\R$, and it is natural to ask which orbits of $\R$ terminate there. For example,
\[
7\longmapsto 8\longmapsto 4\longmapsto 2\longmapsto 1,
\qquad
15\longmapsto 24\longmapsto 12\longmapsto 6\longmapsto 3\longmapsto 4\longmapsto
2\longmapsto 1,
\]
whereas the orbit of $9$ increases before it decreases, namely $$9\mapsto13\mapsto14\mapsto7\mapsto8\mapsto4\mapsto2\mapsto1.$$

A direct computation shows that the orbit of every positive integer $n \le 10^{10}$ reaches $1$: for each odd $n$ it was verified, using a segmented sieve for $\sigma$, that some term of its orbit is smaller than $n$, which suffices since all smaller integers have been treated already. This suggests the following conjecture.

\begin{conjecture}\label{conj}
For every positive integer $n$ there exists $k \in \N$ such that $\R^{k}(n)=1$.
\end{conjecture}

The map $\R$ has appeared before in an empirical form: the number of steps required by $n$ to reach $1$ under $\R$ is sequence A259846 in the On-Line Encyclopedia of Integer Sequences \cite{OEIS}, where the orbit of $9$ displayed above occurs as an example. Iterations of $\sigma$ and of $\varphi$ are discussed by Guy \cite[Section~B41]{Guy2004}, and the normal behavior of such iterates was studied by Erd\H{o}s, Granville, Pomerance and Spiro \cite{ErdosGranvillePomeranceSpiro1990}.

The map $\R$ should be distinguished from the iteration of $n \mapsto \sigma(n)-n$, whose orbits are the aliquot sequences: Catalan and Dickson conjectured that every such sequence is bounded, whereas Guy and Selfridge \cite{GuySelfridge1975} conjectured that for asymptotically all even $n$ it is unbounded; see \cite[Section~B6]{Guy2004}. For $\R$ the arithmetic step is followed by the removal of a power of $2$, and it is this removal which produces the contraction of Corollary~\ref{cor:sqfree}; for the aliquot sequences there is no analogue of it.

The $2$-adic valuation of $\sigma(n)$, which governs the number of halving steps, has itself been studied. Amdeberhan, Moll, Sharma and Villamizar \cite{AMSV2021} evaluated $v_{2}(\sigma(p^{\alpha}))$ and deduced that $v_{2}(\sigma(n)) \le \lceil\log_{2}n\rceil$, with equality precisely when $n$ is a product of distinct Mersenne primes; their bound for odd primes, which was conditional on a technical hypothesis, was made unconditional by Zhao and Chen \cite{ZhaoChen2025}. These are statements about a single value of $\sigma$, and by themselves they say nothing about orbits. We are not aware of any study of the orbit structure of $\R$ itself.

The interest of Conjecture~\ref{conj} lies in its relation to a much older question, namely the existence of odd perfect numbers. Recall that a positive integer $N$ is said to be \emph{perfect} if $\sigma(N)=2N$. If $N$ is an odd perfect number, then $N \mapsto 2N \mapsto N \mapsto 2N \mapsto \cdots$ is a cycle of $\R$ which never reaches $1$. More generally, $N$ is said to be an \emph{odd $2^{k}$-perfect number} if $N$ is odd and $\sigma(N)=2^{k}N$ for some $k \ge 1$, and any such $N$ lies on a cycle of $\R$ of length $k+1$. Consequently, Conjecture~\ref{conj} implies that there is no odd $2^{k}$-perfect number for any $k \ge 1$, and in particular that there is no odd perfect number.

Euler \cite{Euler1849} showed that an odd perfect number must have the form $N=q^{\alpha}s^{2}$, where $q$ is a prime, $\gcd(q,s)=1$ and $q \equiv \alpha \equiv 1 \pmod 4$. Touchard \cite{Touchard1953} showed that $N \equiv 1 \pmod{12}$ or $N \equiv 9 \pmod{36}$, with simplified proofs in \cite{Raghavachari1966, Satyanarayana1959}, and Ochem and Rao \cite{OchemRao2012} showed that $N>10^{1500}$. For further background, see Guy \cite[Section~B1]{Guy2004}. The range $10^{10}$ of the computation above lies very far below that bound, so what it tests is only the orbit structure of $\R$ in the small.

Conjecture~\ref{conj} is strictly stronger than the nonexistence of odd perfect numbers, since it excludes longer cycles and unbounded orbits as well. Our aim is to obtain quantitative information about the orbits of $\R$ and to locate the difficulty which remains.

In Section~\ref{sec:strip} the study of $\R$ is reduced to that of a map $T$ induced by $\R$ on the odd integers, and a pointwise descent estimate for $T$ is proved (Theorem~\ref{thm:descent}). It gives a contraction by a factor of at most $(2/3)^{\omega(n)}$ for odd squarefree $n$, reduces Conjecture~\ref{conj} to the odd integers divisible by the square of a prime (Theorem~\ref{thm:reduction}), and restricts the cycles of $\R$ (Theorem~\ref{thm:cycles}).

Section~\ref{sec:density} contains the main result of the paper, Theorem~\ref{thm:density-descent-strong}: for every fixed $\delta$ with $0<\delta<1$, almost every odd integer satisfies
\[
T(n)\le n(\log n)^{-c'(1-\delta)}, \qquad c'=2\log 2,
\]
where ``almost every'' is meant in the sense of relative density among the odd integers. The arithmetic input is a lower bound for the strongly additive function $\sum_{p \mid n}v_{2}(p+1)$, obtained from the Tur\'an--Kubilius inequality. The corresponding question for Euler's function was settled by Banks, Luca and Shparlinski \cite[Theorem~8]{BanksLucaShparlinski2006}: $v_{2}(\varphi(n))$ has normal order $2\log\log n$. Proposition~\ref{prop:v2normal} is the analogue which $T$ requires, and the argument is theirs; the two differ in that every prime divisor of $n$ contributes to $v_{2}(\varphi(n))$, whereas only those dividing $n$ to an odd exponent contribute to $v_{2}(\sigma(n))$. That parity is what must be controlled, and it is also what obstructs iteration.

\section{The Induced Map and a Descent Estimate}\label{sec:strip}

For a positive integer $n$, denote by $\odd(n):=n/2^{v_{2}(n)}$ the odd part of $n$, where $v_{p}(n)$ is the exponent of $p$ in $n$, that is, the $p$-adic valuation of $n$. As usual, $\omega(n)$ denotes the number of distinct prime divisors of $n$. If $n=2^{a}u$ with $u$ odd, then $\R^{a}(n)=u$. Thus the behavior of $\R$ reduces to that of the induced map
\[
T(n):=\odd\big(\sigma(n)\big)
\]
on the set of odd positive integers. We make this precise as follows.

\begin{lemma}\label{lem:reduce}
Let $n$ be an odd positive integer and put $v:=v_{2}(\sigma(n))$. Then $\R^{v+1}(n)=T(n)$. Consequently, the odd terms occurring in the $\R$-orbit of $n$ are exactly $n, T(n), T^{2}(n), \dots$, and the $\R$-orbit of an arbitrary positive integer reaches $1$ if and only if the $T$-orbit of its odd part reaches $1$.
\end{lemma}

\begin{proof}
Since $n$ is odd, we have $\R(n)=\sigma(n)=2^{v}T(n)$, where $T(n)$ is odd. The next $v$ applications of $\R$ remove the factor $2^{v}$, and hence $\R^{v+1}(n)=T(n)$. Every term of the orbit which lies strictly between two consecutive odd terms is even. Consequently the odd terms of the orbit are precisely $n, T(n), T^{2}(n), \dots$, and the last assertion follows, since $T(1)=1$.
\end{proof}

The $T$-orbit of an odd integer $n$ reaches $1$ in a single step precisely when $\sigma(n)$ is a power of $2$. It turns out that this condition is quite rigid, and the following characterization of Amdeberhan, Moll, Sharma and Villamizar \cite[Theorem~1.3]{AMSV2021}, which they obtained as the equality case of the bound $v_{2}(\sigma(n)) \le \lceil\log_{2}n\rceil$, describes it completely.

\begin{lemma}\cite[Theorem~1.3]{AMSV2021} \label{lem:pow2}
Let $p$ be an odd prime and let $\alpha \ge 1$ be an integer. Then $\sigma(p^{\alpha})$ is a power of $2$ if and only if $\alpha=1$ and $p$ is a Mersenne prime, that is, $p+1$ is itself a power of $2$. Consequently, for an odd positive integer $n$, we have $T(n)=1$ if and only if $n$ is a product of pairwise distinct Mersenne primes.
\end{lemma}

It follows that the squarefree products of Mersenne primes, namely $3, 7, 21, 31, 93, 217, 651, \dots$, are exactly the odd integers whose $T$-orbit reaches $1$ in a single step; whether this family is infinite is the classical open problem of whether there are infinitely many Mersenne primes.

The heuristic which underlies Conjecture~\ref{conj} is the following. Once the powers of $2$ have been removed, that is, once $\sigma$ has been replaced by $T$, the expanding branch of $\R$ grows only slowly. Indeed, for an odd integer $n$, the quantity $\sigma(n)/n$ is bounded above by $\prod_{p \mid n} p/(p-1)$, which grows very slowly with $n$, whereas the number of factors of $2$ available for cancellation in $\sigma(n)$ grows with the number of distinct prime divisors of $n$. We shall now make this observation precise.

For an odd positive integer $n$, let
\[
\wo(n):=\#\{p : v_{p}(n) \text{ is odd}\}
\qquad \text{and} \qquad
E(n):=\prod_{\substack{p \mid n \\ v_{p}(n) \text{ even}}} \frac{p}{p-1},
\]
with the convention that $E(n)=1$ if no prime divides $n$ to an even power. Note that $\wo(n)=0$ if and only if $n$ is a perfect square.

\begin{theorem}\label{thm:descent}
Let $n>1$ be an odd integer. Then $2^{\wo(n)}$ divides $\sigma(n)$, and moreover
\begin{equation}\label{eq:descent}
\frac{T(n)}{n}<\left(\frac{3}{4}\right)^{\wo(n)}E(n).
\end{equation}
\end{theorem}

\begin{proof}
Let $p$ be an odd prime and let $\alpha \ge 1$ be an integer. The sum $\sigma(p^{\alpha})=\sum_{j=0}^{\alpha} p^{j}$ consists of $\alpha+1$ odd terms, and is therefore even precisely when $\alpha+1$ is even, that is, precisely when $\alpha$ is odd. Since $n$ is odd and $\sigma$ is multiplicative, $\sigma(n)=\prod_{p \mid n}\sigma(p^{v_{p}(n)})$ is a product of $\omega(n)$ factors, of which exactly those with $v_{p}(n)$ odd are even. There are $\wo(n)$ such factors, and hence $2^{\wo(n)}$ divides $\sigma(n)$. Since $T(n)$ is the odd part of $\sigma(n)$, we obtain
\begin{equation}\label{eq:half}
T(n) \le \frac{\sigma(n)}{2^{\wo(n)}} .
\end{equation}

Next, for every prime $p$ and every $\alpha \ge 1$, summation of a geometric series gives
\[
\frac{\sigma(p^{\alpha})}{p^{\alpha}}=\sum_{j=0}^{\alpha}p^{-j}<\sum_{j\ge0}p^{-j}=\frac{p}{p-1},
\]
and therefore $\sigma(n)/n<\prod_{p \mid n}p/(p-1)$, again by multiplicativity. Combining this with \eqref{eq:half}, and separating the product according to the parity of $v_{p}(n)$, we readily see that
\[
\frac{T(n)}{n}<2^{-\wo(n)}\prod_{p \mid n}\frac{p}{p-1}
=\left(\prod_{\substack{p \mid n \\ v_{p}(n) \text{ odd}}}\frac{p}{2(p-1)}\right)E(n),
\]
where the first product has exactly $\wo(n)$ factors. Now the function $p \mapsto p/(2(p-1))=\tfrac12\big(1+1/(p-1)\big)$ is decreasing in $p$, and every prime dividing $n$ is at least $3$. Hence each factor of that product is at most $3/4$, and this yields \eqref{eq:descent}.
\end{proof}

\begin{corollary}\label{cor:sqfree}
Let $n>1$ be an odd squarefree integer with $m$ distinct prime factors. Then $T(n) \le (2/3)^{m} n < n$.
\end{corollary}

\begin{proof}
Write $n=p_{1}\cdots p_{m}$, where $p_{1},\dots,p_{m}$ are distinct odd primes; since $n>1$, we have $m \ge 1$. By multiplicativity, $\sigma(n)=\prod_{j=1}^{m}(p_{j}+1)$, and each $p_{j}+1$ is even. Hence $2^{m}$ divides $\sigma(n)$, so that $T(n) \le \sigma(n)/2^{m}$, and therefore
\[
T(n) \le \prod_{j=1}^{m}\frac{p_{j}+1}{2}
=n\prod_{j=1}^{m}\frac{p_{j}+1}{2p_{j}}
\le \left(\frac{2}{3}\right)^{m}n,
\]
where the last inequality holds because $(p+1)/(2p)=\tfrac12(1+1/p)$ is decreasing in $p$, and is hence at most $(3+1)/(2\cdot3)=2/3$ for every prime $p \ge 3$. Since $m \ge 1$, we obtain $T(n) \le \tfrac23 n<n$.
\end{proof}

\begin{corollary}\label{cor:criterion}
Let $n>1$ be odd. If $(4/3)^{\wo(n)} \ge E(n)$, then $T(n)<n$. This holds automatically when $n$ is squarefree, so the criterion can fail only when $n$ is divisible by the square of a prime.
\end{corollary}

\begin{proof}
The implication is immediate from \eqref{eq:descent}. If $n$ is squarefree then $E(n)=1$ and $\wo(n)=\omega(n) \ge 1$, so that $(4/3)^{\wo(n)} \ge 4/3>1=E(n)$.
\end{proof}

The smallest instances of failure are $n=9$ and $n=25$, for which $T(9)=13$ and $T(25)=31$.

The following theorem shows that the squarefree case accounts for the whole of Conjecture~\ref{conj}.

\begin{theorem}\label{thm:reduction}
Conjecture~\ref{conj} holds for every positive integer if and only if it holds for every odd integer which is divisible by the square of a prime.
\end{theorem}

\begin{proof}
The necessity of the stated condition is trivial, and so we prove only its sufficiency. Assume, then, that the $\R$-orbit of every odd integer divisible by the square of a prime reaches $1$. We shall show by strong induction on $n$ that the $\R$-orbit of every positive integer $n$ reaches $1$.

The case $n=1$ is immediate. Let $n>1$, and assume that the orbit of every positive integer smaller than $n$ reaches $1$. Write $n=2^{a}u$ with $u$ odd. If $a \ge 1$, then $\R^{a}(n)=u<n$, and the induction hypothesis applied to $u$ finishes the argument. We may therefore assume that $n$ is odd. If $n$ is not squarefree, then $n$ is an odd integer divisible by the square of a prime, and its orbit reaches $1$ by hypothesis; there is nothing more to prove in this case. Finally, suppose that $n$ is squarefree. Then Lemma~\ref{lem:reduce} and Corollary~\ref{cor:sqfree} give $\R^{k}(n)=T(n)<n$ for $k=v_{2}(\sigma(n))+1$, and the induction hypothesis, applied to the integer $\R^{k}(n)<n$, shows that the orbit of $\R^{k}(n)$, and hence that of $n$, reaches $1$.
\end{proof}

The strictness of the inequality $\R^{k}(n)<n$ is essential here, since $\R^{k}(n) \le n$ holds trivially for every element of a cycle.

The class left open is not a thin one: the odd squarefree integers have relative density $8/\pi^{2}$ among the odd integers, so that the integers left open by Theorem~\ref{thm:reduction} have relative density $1-8/\pi^{2}=0.1894\dots$ there. The reduction removes rather more than four fifths of the odd integers, but it does not reduce the problem to a set of relative density zero.

Corollary~\ref{cor:sqfree} also yields an explicit family of integers for which the conjecture can be verified directly, provided that squarefreeness persists along the orbit. Call an odd positive integer $n$ \emph{tame} if every term of the sequence $n, T(n), T^{2}(n), \dots$ is squarefree. If $n$ is tame, then Corollary~\ref{cor:sqfree} gives $T^{i+1}(n) \le \tfrac23 T^{i}(n)$ as long as $T^{i}(n)>1$, and hence $T^{j}(n)=1$ for some $j \le \log_{3/2}n$. By Lemma~\ref{lem:reduce}, the orbit of $n$ under $\R$ then reaches $1$ in $j+\sum_{i<j}v_{2}(\sigma(T^{i}(n)))$ steps, the extra terms counting the intermediate halvings. Tameness can therefore be tested using at most $O(\log n)$ factorizations.

An explicit family of tame integers is obtained from chains of primes. Let $P_{1}$ denote the set of Mersenne primes, and for $i \ge 2$ put $P_{i}:=\{p \text{ prime} : \odd(p+1) \in P_{i-1}\}$, so that, for instance, $P_{2}=\{5,11,13,23,47,61,\dots\}$ and $P_{3}=\{19,43,79,103,367,\dots\}$. If $p \in P_{i}$, then $T(p)=\odd(p+1) \in P_{i-1}$, and hence the $T$-orbit of $p$ is a chain of primes which terminates at $1$, by Lemma~\ref{lem:pow2}. Since every prime is squarefree, every element of $\bigcup_{i}P_{i}$ is tame.

\begin{remark}
Tameness does not extend automatically from primes to their products. Let $n=p_{1}\cdots p_{m}$, where $p_{j}+1=2^{r_{j}}q_{j}$ with $q_{j}$ odd; then multiplicativity gives $T(n)=q_{1}\cdots q_{m}$. The $q_{j}$ need not be distinct, so that $T(n)$ need not be squarefree even when $n$ is. For example, $5$ and $11$ both lie in $P_{2}$, while $\sigma(55)=6\cdot12=72=2^{3}\cdot3^{2}$, so that $T(55)=9$ is not squarefree. The orbit of $55$ nonetheless reaches $1$, since $T(9)=13 \in P_{2}$, but this has to be checked separately. What holds unconditionally is that if $q_{1},\dots,q_{m}$ are pairwise distinct and $q_{1}\cdots q_{m}$ is tame, then $n$ is tame as well.
\end{remark}

By a \emph{cycle} of $\R$ we mean the orbit $\{n,\R(n),\R^{2}(n),\dots\}$ of a positive integer $n$ for which $\R^{k}(n)=n$ for some $k \ge 1$; such a set is finite and is permuted cyclically by $\R$. Recall now that an orbit of $\R$ can fail to reach $1$ in one of two ways only, namely by entering a cycle, or by increasing without bound. The following theorem concerns the former possibility.

\begin{theorem}\label{thm:cycles}
Let $C \ne \{1\}$ be a cycle of $\R$, and let $m_{1},\dots,m_{\ell}$ denote the odd elements of $C$, indexed cyclically so that $T(m_{i})=m_{i+1}$ for each $i$, where the indices are read modulo $\ell$. Put $v_{i}:=v_{2}(\sigma(m_{i}))$, let $V:=v_{1}+\cdots+v_{\ell}$, and let $I(n):=\sigma(n)/n$ denote the abundancy index of $n$. Then $\ell \ge 1$, and the following statements hold.
\begin{enumerate}[label=\textup{(\alph*)}, leftmargin=2.4em]
\item $\displaystyle \prod_{i=1}^{\ell} I(m_{i})=2^{V}$.
\item $\ell=1$ if and only if $m_{1}$ is an odd $2^{v_{1}}$-perfect number. The case $v_{1}=1$ is precisely that of an odd perfect number.
\item The smallest of the integers $m_{1},\dots,m_{\ell}$ is not squarefree.
\item If none of $m_{1},\dots,m_{\ell}$ is a perfect square, then $I(m_{i}) \ge 2$ for at least one index $i$.
\end{enumerate}
\end{theorem}

\begin{proof}
Every even element of $C$ is mapped by $\R$ to a strictly smaller integer, and hence $C$ cannot consist of even elements alone. Thus $\ell \ge 1$. By Lemma~\ref{lem:reduce}, the odd elements of $C$ are permuted cyclically by $T$, and this is the indexing which we have fixed.

For (a), the definition of $T$ gives $\sigma(m_{i})=2^{v_{i}}m_{i+1}$ for each $i$. Multiplying these $\ell$ equations, we obtain
\[
\prod_{i=1}^{\ell}\sigma(m_{i})=2^{V}\prod_{i=1}^{\ell}m_{i+1}=2^{V}\prod_{i=1}^{\ell}m_{i},
\]
where the last equality holds because the indices are read modulo $\ell$, so that $m_{2},\dots,m_{\ell},m_{1}$ is a permutation of $m_{1},\dots,m_{\ell}$. Division by $\prod_{i}m_{i}$, which is non-zero, gives $\prod_{i}I(m_{i})=2^{V}$.

Statement (b) is immediate from the definitions. Indeed, $\ell=1$ means that $T(m_{1})=m_{1}$, that is, $\sigma(m_{1})=2^{v_{1}}m_{1}$ with $m_{1}$ odd, and this says precisely that $m_{1}$ is an odd $2^{v_{1}}$-perfect number. The case $v_{1}=1$ is that of an odd perfect number.

To see (c), let $m_{i_{0}}$ be the smallest of the integers $m_{1},\dots,m_{\ell}$. Its image $T(m_{i_{0}})=m_{i_{0}+1}$ again belongs to $\{m_{1},\dots,m_{\ell}\}$, and hence $T(m_{i_{0}}) \ge m_{i_{0}}$, by minimality. Moreover $m_{i_{0}}>1$, since $C \ne \{1\}$. Were $m_{i_{0}}$ squarefree, Corollary~\ref{cor:sqfree} would give $T(m_{i_{0}})<m_{i_{0}}$. Hence $m_{i_{0}}$ is not squarefree.

Finally, we have seen in the proof of Theorem~\ref{thm:descent} that $\sigma(n)$ is odd exactly when $\wo(n)=0$, that is, exactly when the odd integer $n$ is a perfect square; this is also \cite[Corollary~3.3]{AMSV2021}. Thus, if none of $m_{1},\dots,m_{\ell}$ is a perfect square, then $v_{i} \ge 1$ for every $i$, and hence $V \ge \ell$. By (a), $\prod_{i}I(m_{i})=2^{V} \ge 2^{\ell}$, and so the $\ell$ factors $I(m_{i})$ cannot all be smaller than $2$. This is (d).
\end{proof}

Statement (b) makes precise the connection with odd perfect numbers described in Section~\ref{sec:intro}, and statement (c) forces any cyclic counterexample into the class which Theorem~\ref{thm:reduction} leaves open, in agreement with Euler's form $N=q^{\alpha}s^{2}$. Thus fixed points and longer cycles present the same difficulty.

\begin{corollary}\label{cor:opn-constraint}
Let $N$ be an odd $2^{k}$-perfect number. Then $E(N)>(4/3)^{\wo(N)}$. If $k=1$, then $\wo(N)=1$ and
\[
\prod_{p \mid s}\frac{p}{p-1}=E(N)>\frac43,
\]
where $N=q^{\alpha}s^{2}$ is Euler's form. This is weaker than the bound $E(N)>2(1-1/q) \ge 8/5$ which follows from $\sigma(N)=2N$ alone.
\end{corollary}

\begin{proof}
Since $N$ is odd and $\sigma(N)=2^{k}N$, we have $T(N)=\odd(2^{k}N)=N$, and \eqref{eq:descent} gives the first inequality. Let $k=1$. By Euler's theorem, quoted in Section~\ref{sec:intro}, $N=q^{\alpha}s^{2}$ with $q \equiv \alpha \equiv 1 \pmod 4$ and $\gcd(q,s)=1$; thus $\alpha$ is odd and every prime dividing $s$ divides $N$ to an even exponent, so that $\wo(N)=1$ and $E(N)=\prod_{p \mid s}p/(p-1)$. For the last assertion, multiplicativity gives $2=I(q^{\alpha})I(s^{2})<\frac{q}{q-1}E(N)$, by the bound $I(m)<\prod_{p \mid m}p/(p-1)$ of the proof of Theorem~\ref{thm:descent}; hence $E(N)>2(1-1/q) \ge 8/5$, since $q \ge 5$.
\end{proof}

Corollary~\ref{cor:opn-constraint} measures how far \eqref{eq:descent} is from the odd perfect number problem: it is weaker than two lines of the definition of perfection. It is the class of integers of Theorem~\ref{thm:reduction}, and not the estimate itself, which carries the difficulty.
Finally, the heuristic advantage of $\R$ over the classical Collatz map can be stated in one line: by Corollary~\ref{cor:sqfree} a squarefree odd $n$ contracts by a factor of at most $(2/3)^{\omega(n)}$, which tends to $0$ since $\omega(n)$ has normal order $\log\log n$ \cite{HardyRamanujan1917}, whereas for the Collatz map the corresponding factor is the fixed constant $3/4$.
\section{An Almost-All Descent Estimate}\label{sec:density}

The map $T$ contracts almost every odd integer by a power of $\log n$ in a single step. We shall prove this in the sharpest form which the present method gives, namely with the exponent $2\log2$ (Theorem~\ref{thm:density-descent-strong}). The estimate concerns one step only. To iterate it, one needs information about the distribution of $T(n)$, and not merely about that of $n$. We shall then determine why this information is not available, and what is known unconditionally about the prime factors of $\sigma(n)$.

The map $T$ of Section~\ref{sec:strip} is defined on the odd positive integers, and by Lemma~\ref{lem:reduce} the behavior of $\R$ is determined by that of $T$. With this in view, all density statements in this section are taken relative to the odd integers.

Write $\No:=\{n \in \N : n \text{ odd}\}$. For $A\subseteq\N$, the \emph{natural density} of $A$ is
\[
d(A):=\lim_{x\to\infty}\frac{\#(A\cap[1,x])}{x},
\]
when the limit exists, and for $A \subseteq \No$ the \emph{relative density of $A$ among the odd integers} is
\[
\dodd(A):=\lim_{x\to\infty}\frac{\#(A\cap[1,x])}{\#(\No\cap[1,x])},
\]
again when the limit exists. Replacing the limit by $\limsup$ or $\liminf$ gives the upper and the lower relative densities $\overline{\dodd}$ and $\underline{\dodd}$. We say that a property holds for \emph{almost all odd $n$} if the set of odd $n$ for which it fails has relative density $0$ among the odd integers. Since $\#(\No\cap[1,x])=\lceil x/2\rceil\sim x/2$, every set of natural density zero has relative density zero after restriction to $\No$:
\begin{equation}\label{eq:transfer}
d(A)=0 \quad\Longrightarrow\quad \dodd(A\cap\No)=0 .
\end{equation}
This allows us to transfer to the odd integers those statements which the literature proves for all $n$.

We write $p_{i}$ for the $i$-th prime, so that $p_{1}=2$ and $p_{2}=3$. If $\mathcal{P}$ is a set of primes, its \emph{relative density among primes} is $\lim_{x\to\infty}\pi(x;\mathcal{P})/\pi(x)$, where $\pi(x;\mathcal{P}):=\#\{p\le x: p\in\mathcal{P}\}$ and $\pi(x)$ is the usual prime-counting function. By Dirichlet's theorem on primes in arithmetic progressions, every residue class $a\bmod m$ with $\gcd(a,m)=1$ contains infinitely many primes. The quantitative facts which we shall use are the prime number theorem for arithmetic progressions, by which such a class has relative density $1/\varphi(m)$ among primes, where $\varphi$ is Euler's totient function, and the estimate
\begin{equation}\label{eq:mertens-ap}
\sum_{\substack{q\le x \\ q\equiv a \!\!\pmod m}}\frac1q
=\frac{\log\log x}{\varphi(m)}+O_{m}(1),
\end{equation}
where the sum is over primes $q$; in particular, this sum diverges. Both of these may be found in \cite[Chapter~20]{Davenport2000}, the second one by partial summation from the first.

For the ordinary Collatz map $C$, writing $C_{\min}(n):=\inf_{k}C^{k}(n)$, Terras \cite[Theorem~1.17]{Terras1976} showed that $C_{\min}(n)<n$ for almost all $n$, and Korec \cite[Theorem~1]{Korec1994} that $C_{\min}(n)\le n^{\theta}$ for every fixed $\theta>\log3/\log4$, both in natural density. Those results concern an entire orbit and ours a single step, so the two are not comparable; they differ in that the Collatz contraction factor at an odd step is the constant $3/4$, whereas here the number of factors of $2$ removed grows with $n$.

By Theorem~\ref{thm:descent} at least one factor of $2$ is available at each prime dividing $n$ to an odd exponent. In fact $v_{2}(p+1)$ are available at such a prime $p$, and the average of $v_{2}(p+1)$ over the primes is $2$; this produces the exponent $2\log2$. The tool is the Tur\'an--Kubilius inequality, and the input is the exact evaluation of $v_{2}(\sigma(n))$, which we now quote.

\begin{theorem} \cite[Theorem~3.2]{AMSV2021}\label{thm:AMSV}
Let $p$ be an odd prime and let $\alpha \ge 1$ be an integer. Then
\[
v_{2}\big(\sigma(p^{\alpha})\big)=
\begin{cases}
0, & \text{if $\alpha$ is even},\\[2pt]
v_{2}(\alpha+1)+v_{2}(p+1)-1, & \text{if $\alpha$ is odd}.
\end{cases}
\]
\end{theorem}

Since $\sigma$ is multiplicative, Theorem~\ref{thm:AMSV} gives, for every odd positive integer $n$,
\begin{equation}\label{eq:v2exact}
v_{2}(\sigma(n))=\sum_{\substack{p \mid n \\ v_{p}(n) \text{ odd}}}\Big(v_{2}\big(v_{p}(n)+1\big)+v_{2}(p+1)-1\Big);
\end{equation}
this identity is \cite[Theorem~1.1]{AMSV2021}, written there in the equivalent form $\sum v_{2}\big((v_{p}(n)+1)(p+1)/2\big)$.
We shall use \eqref{eq:v2exact} only through the following consequence, in which the terms $v_{2}(v_{p}(n)+1)-1$ are discarded.

\begin{corollary}\label{cor:v2lower}
For every odd positive integer $n$ we have $v_{2}(\sigma(n))=F(n)+D(n)$, where
\begin{equation}\label{eq:v2lower}
F(n):=\sum_{\substack{p \mid n \\ v_{p}(n) \text{ odd}}} v_{2}(p+1),
\qquad
D(n):=\sum_{\substack{p \mid n \\ v_{p}(n) \text{ odd}}}\Big(v_{2}\big(v_{p}(n)+1\big)-1\Big).
\end{equation}
Moreover $D(n) \ge 0$, and hence $v_{2}(\sigma(n)) \ge F(n)$.
\end{corollary}

\begin{proof}
The identity is \eqref{eq:v2exact}. If $v_{p}(n)$ is odd, then $v_{p}(n)+1$ is even, so that $v_{2}(v_{p}(n)+1) \ge 1$ and the corresponding term of $D(n)$ is non-negative. Hence $D(n) \ge 0$.
\end{proof}

It may be noted that $D(n)>0$ only when some prime divides $n$ to an exponent $\alpha$ with $\alpha \equiv 3 \pmod 4$; in particular, $D(n)=0$ unless $n$ is divisible by the cube of a prime. The divisibility $2^{\wo(n)} \mid \sigma(n)$ of Theorem~\ref{thm:descent} is also immediate from \eqref{eq:v2exact}, since every term of that sum is at least $1$.

The quantity $F$ is not additive, since it depends on the parity of the exponents in $n$. We therefore compare it with the function
\[
G(n):=\sum_{p \mid n} v_{2}(p+1),
\]
which depends only on the set of primes dividing $n$ and is therefore strongly additive. It is easier to handle than $F$, and we bound the difference in mean.

\begin{lemma}\label{lem:v2correction}
For every $x \ge 1$,
\[
\sum_{\substack{n \le x \\ n \text{ odd}}}\big(G(n)-F(n)\big) \le C_{1}\,x, \qquad
C_{1}:=\sum_{p>2}\frac{v_{2}(p+1)}{p^{2}-1}<\infty ,
\]
where the sum defining $C_{1}$ is over all odd primes $p$.
\end{lemma}

\begin{proof}
By the definitions of $F$ and $G$, we have $G(n)-F(n)=\sum v_{2}(p+1)$, where the sum is over those primes $p$ which divide $n$ to an even exponent. Such a prime is odd, since $n$ is odd, and $v_{p}(n)=2a$ for some integer $a \ge 1$. Writing this contribution as an indicator and summing over $n$, we obtain
\[
\sum_{\substack{n \le x \\ n \text{ odd}}}\big(G(n)-F(n)\big)
=\sum_{p>2}v_{2}(p+1)\sum_{a \ge 1}\#\{n \le x : n \text{ odd},\ v_{p}(n)=2a\},
\]
where the interchange of the summations is legitimate, since all the terms are non-negative. If $v_{p}(n)=2a$, then $p^{2a}$ divides $n$, and hence the innermost count is at most $x/p^{2a}$. Summation of the resulting geometric series gives $\sum_{a \ge 1}x/p^{2a}=x/(p^{2}-1)$, and therefore the left-hand side is at most $C_{1}x$. Finally, $C_{1}$ is finite. Indeed, $2^{v_{2}(p+1)} \le p+1$, so that $v_{2}(p+1) \le \log(p+1)/\log 2$, and $\sum_{p}\log(p+1)/(p^{2}-1)$ converges.
\end{proof}

\begin{corollary}\label{cor:v2correction}
For every $\varepsilon>0$ there is an integer $K_{1}=K_{1}(\varepsilon) \ge 1$ such that
\[
\#\{n \le x : n \text{ odd},\ G(n)-F(n) \ge K_{1}\} \le \varepsilon\,\#(\No\cap[1,x])
\qquad \text{for every } x \ge 1.
\]
\end{corollary}

\begin{proof}
Since $G-F$ takes non-negative values, discarding from the sum in Lemma~\ref{lem:v2correction} all the terms with $G(n)-F(n)<K_{1}$ only decreases it, and each remaining term is at least $K_{1}$. This is Markov's inequality in the present discrete setting, and it gives $\#\{n \le x : n \text{ odd},\ G(n)-F(n) \ge K_{1}\} \le C_{1}x/K_{1}$. Since $\#(\No\cap[1,x])=\lceil x/2\rceil \ge x/2$, we have $x \le 2\,\#(\No\cap[1,x])$, and it remains to choose $K_{1}$ with $2C_{1}/K_{1}<\varepsilon$.
\end{proof}

The present argument needs a bound on the abundancy index $I(n)=\sigma(n)/n$, and such a bound is available in mean.

\begin{lemma}\label{lem:abundancy-mean}
For every $x \ge 1$ we have $\sum_{n \le x} \sigma(n)/n \le (\pi^{2}/6)\,x$. Consequently, for every $\varepsilon>0$ there is a constant $C_{2}=C_{2}(\varepsilon)\ge1$ such that
\[
\#\{n \le x : n \text{ odd},\ \sigma(n)/n \ge C_{2}\} \le \varepsilon\,\#(\No\cap[1,x])
\qquad \text{for every } x \ge 1.
\]
\end{lemma}

\begin{proof}
Since $\sigma(n)/n=\sum_{d \mid n}1/d$, interchanging the order of summation gives
\[
\sum_{n \le x}\frac{\sigma(n)}{n}=\sum_{d \le x}\frac{1}{d}\Big\lfloor\frac{x}{d}\Big\rfloor
\le x\sum_{d \ge 1}\frac{1}{d^{2}}=\frac{\pi^{2}}{6}\,x .
\]
The odd integers $n \le x$ form a subset of $[1,x]$ and $\sigma(n)/n>0$, so the same bound holds for the sum restricted to them. Markov's inequality, together with $x \le 2\,\#(\No\cap[1,x])$, now gives the second assertion with $C_{2}:=\lceil \pi^{2}/(3\varepsilon)\rceil+1$.
\end{proof}

We come to a lower bound for $G$ which holds for almost all $n$. The argument is that of Banks, Luca and Shparlinski \cite[Theorem~8]{BanksLucaShparlinski2006}, who determined the normal order of $v_{q}(\varphi(n))$ for every prime $q$; we need only a lower bound, and a truncation makes it elementary. For a fixed integer $J \ge 1$, put
\[
G_{J}(n):=\sum_{p \mid n}\min\{v_{2}(p+1),\,J\}, \qquad \kappa_{J}:=2-2^{1-J},
\]
so that $G_{J}$ is strongly additive, $0 \le G_{J}(p) \le J$ for every prime $p$, and $G_{J}(n) \le G(n)$ for every $n$. The truncation is what makes the following proposition elementary; the constant $\kappa_{J}$ tends to $2$ as $J \to \infty$, and this is the source of the exponent $2\log 2$.

\begin{proposition}\label{prop:v2normal}
Let $J \ge 1$ be fixed. Then
\[
G_{J}(n)\ \ge\ \kappa_{J}\log\log n-(\log\log n)^{3/4}
\]
for almost all $n>e^{e}$, in natural density, and hence, by \eqref{eq:transfer}, for almost all odd $n$.
\end{proposition}

The restriction $n>e^{e}$ serves only to make $\log\log n$ positive. It excludes finitely many integers and therefore affects no density.

\begin{proof}
For a non-negative integer $v$ we have $\min\{v,J\}=\sum_{k=1}^{J}\mathbf{1}_{\{v \ge k\}}$, and $v_{2}(p+1) \ge k$ precisely when $p \equiv -1 \pmod{2^{k}}$. Hence, on putting
\[
A_{J}(x):=\sum_{p \le x}\frac{G_{J}(p)}{p},
\]
we obtain
\[
A_{J}(x)=\sum_{k=1}^{J}\ \sum_{\substack{p \le x \\ p \equiv -1 \!\!\pmod{2^{k}}}}\frac1p
=\sum_{k=1}^{J}\Big(\frac{\log\log x}{2^{k-1}}+O_{k}(1)\Big)
=\kappa_{J}\log\log x+O_{J}(1),
\]
where we have used \eqref{eq:mertens-ap} with $m=2^{k}$ and $a=-1$, which is permissible since $\gcd(-1,2^{k})=1$, together with $\varphi(2^{k})=2^{k-1}$ and $\sum_{k=1}^{J}2^{1-k}=2-2^{1-J}=\kappa_{J}$. Moreover, since $0 \le G_{J}(p) \le J$,
\[
B_{J}(x)^{2}:=\sum_{p \le x}\frac{G_{J}(p)^{2}}{p}\ \le\ J\,A_{J}(x)\ \ll_{J}\ \log\log x .
\]
The mean and the variance occurring in the Tur\'an--Kubilius inequality are sums over prime powers $p^{\nu} \le x$. For a strongly additive function they differ from $A_{J}(x)$ and $B_{J}(x)^{2}$ by $O_{J}(1)$ in each case, since $\sum_{\nu \ge 1}p^{-\nu}=1/(p-1)$ and $0 \le G_{J}(p) \le J$; such differences are absorbed into the error terms above. That inequality (see, for example, Tenenbaum \cite[Chapter~III.3]{Tenenbaum2015}) now gives
\[
\sum_{n \le x}\big|G_{J}(n)-A_{J}(x)\big|^{2}\ \ll\ x\,B_{J}(x)^{2}\ \ll_{J}\ x\log\log x .
\]
Consequently, by Chebyshev's inequality,
\[
\#\Big\{n \le x : \big|G_{J}(n)-A_{J}(x)\big| \ge \tfrac12(\log\log x)^{3/4}\Big\}
\ \ll_{J}\ \frac{x\log\log x}{(\log\log x)^{3/2}}=\frac{x}{(\log\log x)^{1/2}}=o(x).
\]

Let now $n \le x$ be an integer lying neither in $[1,\sqrt{x}\,]$ nor in the exceptional set displayed above; the first of these has $\sqrt{x}=o(x)$ elements, and the second has $o(x)$ elements, so that all $n \le x$ outside a set of $o(x)$ elements are of this kind. For such an $n$ we have $n>\sqrt{x}$, and hence $\log\log n \ge \log\log x-\log 2$. Therefore
\[
G_{J}(n)\ \ge\ A_{J}(x)-\tfrac12(\log\log x)^{3/4}
\ \ge\ \kappa_{J}\log\log n-\tfrac12(\log\log x)^{3/4}-O_{J}(1),
\]
where we have used $\log\log n \le \log\log x$ together with $\kappa_{J}>0$. Since $\log\log x \le \log\log n+\log 2$, the quantity $\tfrac12(\log\log x)^{3/4}+O_{J}(1)$ is at most $(\log\log n)^{3/4}$ for all $n$ larger than a bound depending only on $J$. The integers below that bound are finite in number and affect no density, and the assertion follows.
\end{proof}

\begin{theorem}\label{thm:density-descent-strong}
Let $c':=2\log2$. For every $\delta$ with $0<\delta<1$, the set
\[
G_{\delta}:=\big\{n \in \No,\ n\ge3 : T(n) \le n\,(\log n)^{-c'(1-\delta)}\big\}
\]
has relative density $1$ among the odd positive integers; that is, the inequality $T(n)\le n(\log n)^{-c'(1-\delta)}$ holds for almost all odd $n$.
\end{theorem}

\begin{proof}
Fix $\delta \in (0,1)$ and $\varepsilon>0$. Choose an integer $J$ with $2^{1-J} \le \delta/2$, so that $\kappa_{J} \ge 2-\delta/2$, and let $K_{1}=K_{1}(\varepsilon)$ and $C_{2}=C_{2}(\varepsilon)$ be as in Corollary~\ref{cor:v2correction} and Lemma~\ref{lem:abundancy-mean}. Put
\[
H_{1}:=\{n \in \No : G(n)-F(n)<K_{1}\}, \qquad
H_{2}:=\{n \in \No : \sigma(n)/n<C_{2}\},
\]
and let $A'$ denote the set of odd $n$ for which the conclusion of Proposition~\ref{prop:v2normal} holds. Then $\underline{\dodd}(H_{1}) \ge 1-\varepsilon$ and $\underline{\dodd}(H_{2}) \ge 1-\varepsilon$, while $\dodd(A')=1$, and therefore
\begin{equation}\label{eq:goodset-strong}
\underline{\dodd}\big(A' \cap H_{1} \cap H_{2}\big)\ \ge\ 1-2\varepsilon ,
\end{equation}
since, on dividing by $\#(\No\cap[1,x])$, the two terms subtracted have upper limit at most $\varepsilon$ each while the third tends to $0$.

Let $n \in A' \cap H_{1} \cap H_{2}$ with $n>e^{e}$. Since $G_{J} \le G$, Corollary~\ref{cor:v2lower} and Proposition~\ref{prop:v2normal} give
\[
v_{2}(\sigma(n))\ \ge\ F(n)\ =\ G(n)-\big(G(n)-F(n)\big)\ \ge\ G_{J}(n)-K_{1}
\ \ge\ \kappa_{J}\log\log n-(\log\log n)^{3/4}-K_{1}.
\]
Since $T(n)=\sigma(n)2^{-v_{2}(\sigma(n))}$ and $\sigma(n)<C_{2}n$, it follows that
\[
\frac{T(n)}{n}\ <\ C_{2}\,2^{K_{1}}\ 2^{(\log\log n)^{3/4}}\ 2^{-\kappa_{J}\log\log n} .
\]
Now $2^{-\kappa_{J}\log\log n}=e^{-\kappa_{J}(\log2)\log\log n}=(\log n)^{-\kappa_{J}\log2}$, and $\kappa_{J}\log 2 \ge (2-\delta/2)\log2=c'(1-\delta/4)$. Hence
\[
\frac{T(n)}{n}\ <\ C_{2}\,2^{K_{1}}\ 2^{(\log\log n)^{3/4}}\,(\log n)^{-c'(1-\delta/4)} .
\]
It remains to absorb the first two factors. Since $(\log\log n)^{3/4}=o(\log\log n)$, we have $2^{(\log\log n)^{3/4}}=(\log n)^{o(1)}$, and therefore both $C_{2}2^{K_{1}} \le (\log n)^{c'\delta/8}$ and $2^{(\log\log n)^{3/4}} \le (\log n)^{c'\delta/8}$ hold for all $n$ larger than some $n_{2}$ depending on $\delta$ and $\varepsilon$ alone. For such $n$ we obtain
\[
\frac{T(n)}{n}\ <\ (\log n)^{-c'(1-\delta/4)+c'\delta/4}=(\log n)^{-c'(1-\delta/2)}\ \le\ (\log n)^{-c'(1-\delta)} .
\]
Thus the set in question contains $(A' \cap H_{1} \cap H_{2}) \cap [n_{2},\infty)$. Removal of the finitely many integers below $n_{2}$ alters no relative density, and so \eqref{eq:goodset-strong} bounds its lower relative density below by $1-2\varepsilon$. Here $\varepsilon>0$ was arbitrary and the set does not depend on $\varepsilon$, so that its relative density is $1$.
\end{proof}

The bound $v_{2}(\sigma(n)) \ge F(n)$ of Corollary~\ref{cor:v2lower} loses very little. The loss is exactly $D(n)$, and a direct computation over all odd integers $3 \le n \le 10^{7}$ gives $0.0359$ for its mean and $3$ for its maximum. Retaining $D(n)$ would therefore not improve the exponent $c'$. It may also be noted that the same argument, with the divisibility $2^{\wo(n)} \mid \sigma(n)$ of Theorem~\ref{thm:descent} in place of Theorem~\ref{thm:AMSV} and the Hardy--Ramanujan theorem in place of Proposition~\ref{prop:v2normal}, yields only the weaker exponent $\log(4/3)$.

We turn now to the obstruction to iterating this estimate. For $A\subseteq\N$, the \emph{logarithmic density} of $A$ is
\[
\lim_{x\to\infty}\frac{\displaystyle\sum_{\substack{n\in A\\ n\le x}}\frac1n}
{\displaystyle\sum_{n\le x}\frac1n},
\]
when the limit exists. A set of natural density $1$ has logarithmic density $1$, but the converse can fail, so that logarithmic density is the weaker of the two notions. Tao's almost-all theorem \cite[Theorem~1.3]{Tao2022} states that for every function $f:\N\to\mathbb{R}$ with $\lim_{N\to\infty}f(N)=+\infty$, one has $C_{\min}(N)<f(N)$ for almost all $N$, in the sense of logarithmic density. That is not the notion used in Theorem~\ref{thm:density-descent-strong}, where the densities are relative densities among the odd integers.

We can now say what an iteration of Theorem~\ref{thm:density-descent-strong} along the orbit $n, T(n), T^{2}(n),\dots$ would require, and which implication fails. Theorem~\ref{thm:density-descent-strong} asserts that $\dodd(G_{\delta})=1$, whereas a second application of the same estimate, now at the point $T(n)$, requires $\dodd\big(T^{-1}(G_{\delta})\big)=1$. The first does not imply the second. No asymptotic preservation of relative density under $T$ is known, so $T$ may carry a set of full relative density into one of relative density zero, where neither Proposition~\ref{prop:v2normal} nor Theorem~\ref{thm:density-descent-strong} says anything. Nothing proved here controls the distribution of $F$ along the image of $T$, and a larger exponent would not help, since what fails is the transfer through $T$ and not the strength of the single-step estimate. This is the analogue of the obstruction in Tao's treatment of the Collatz map, which he overcomes by a $3$-adic Fourier-analytic argument \cite{Tao2022}; whether such an argument can be adapted to $T$ remains open.

It remains to say what is known unconditionally about the prime factors of $\sigma(n)$, since these govern $\omega(T(n))$, and to determine which part of the obstruction described above it removes. The mechanism is elementary: if $q \equiv -1 \pmod p$ and $v_{q}(n)$ is odd, then $q+1$, and hence $p$, divides $\sigma(q^{v_{q}(n)})$, and $\sum_{q \equiv -1 (p)}1/q$ diverges by \eqref{eq:mertens-ap}. A far stronger statement is available unconditionally: Erd\H{o}s and Pomerance \cite{ErdosPomerance1985} determined the normal number of prime factors of $\varphi$, and their methods give, as noted by Troupe \cite[\S 5]{Troupe2015}, the corresponding statement for $\sigma$, namely that $\omega(\sigma(n))$ is asymptotically normally distributed with mean $\tfrac12(\log\log n)^{2}$ and standard deviation $\tfrac{1}{\sqrt3}(\log\log n)^{3/2}$. In particular $\omega(\sigma(n))$, and hence $\omega(T(n))$, exceeds $(1-\eta)\tfrac12(\log\log n)^{2}$ for almost all $n$, for every fixed $\eta>0$.

This does not suffice. A second application of Theorem~\ref{thm:density-descent-strong} at the point $T(n)$ requires not $\omega(T(n))$ but $F(T(n))$, the sum of $v_{2}(p+1)$ over the primes dividing $T(n)$ to an \emph{odd} exponent. The result just quoted counts the prime factors of $T(n)$ and says nothing about the parity of their exponents, and it is exactly that parity on which \eqref{eq:v2lower} depends.

A second difficulty is independent of the one just described. The normal-order results quoted above describe $n$ as it ranges over $[1,x]$, whereas a second iteration of Theorem~\ref{thm:density-descent-strong} requires control of $F(T(n))$ and of the abundancy index of $T(n)$ relative to the image of $T$, and hence the distribution of $v_{2}(\sigma(T(n)))$ there. Neither is treated in the classical normal-order literature. Together with the failure described above, these questions constitute the principal remaining obstacle to an unconditional almost-all statement for Conjecture~\ref{conj}, and we leave them open for future investigation.

\bibliographystyle{plain}
\bibliography{refs}

\end{document}